\documentclass{amsart}
\title[An Approach to Hopf Algebras II]{An Approach to Hopf Algebras via Frobenius Coordinates II}
\author{Lars Kadison and A.A. Stolin}
\address{Chalmers University of Technology/G{\" o}teborg
University \\
S-412 96 G{\" o}teborg, Sweden}
\email{kadison@math.ntnu.no \\
astolin@math.chalmers.se}
\thanks{The  first author thanks NorFA in Oslo for financial support.}  
\date{}

\newtheorem{theorem}{Theorem}[section]
\newtheorem{lemma}[theorem]{Lemma}
\newtheorem{proposition}[theorem]{Proposition}

\theoremstyle{definition}
\newtheorem{definition}[theorem]{Definition}

\newtheorem{corollary}[theorem]{Corollary}

\theoremstyle{remark}
\newtheorem{remark}[theorem]{Remark}

\newcommand{\id}{\mbox{id}}

\newcommand{\End}{\mbox{End}}
\newcommand{\Hom}{\mbox{Hom}}

\newcommand{\eps}{\varepsilon}
\newcommand{\es}{\overline{S}}

\newcommand{\flip}{\varsigma}

\newcommand{\0}{_{(0)}}
\newcommand{\1}{_{(1)}}
\newcommand{\2}{_{(2)}}
\newcommand{\3}{_{(3)}}

\newcommand{\I}{_{i(1)}}
\newcommand{\II}{_{i(2)}}
\newcommand{\III}{_{i(3)}}

\input diagrams
\begin{document}
\maketitle
\begin{abstract}
We study a Hopf algebra $H$, which is finitely generated and projective
over a commutative ring $k$, as a $P$-Frobenius algebra.  We define 
modular functions in this setting, and provide a complete proof of
Radford's formula for the fourth power of the antipode, using Frobenius
algebraic techniques.  As further applications, we extend Etingof and
Gelaki's result that a separable and coseparable Hopf algebra has
antipode of order two, the result of Schneider that Hopf subalgebras
are twisted Frobenius extensions, and show that the quantum double
is always a Frobenius algebra.  
\end{abstract}

\section{Introduction}

Perhaps the first beginnings of relating Frobenius algebras to Hopf algebras
was the example
by Berkson \cite{Berk}. He proved that the restricted universal enveloping algebra of
a finite dimensional restricted Lie algebra is a Frobenius
algebra. Together with the well-known Frobenius
algebra examples 
of finite group algebras, this raised the question if 
a finitely generated, projective Hopf $k$-algebra $H$ is Frobenius. This  
was established 
by Larson and Sweedler \cite{LS} for $k$ a principal ideal domain
and was generalized by Pareigis \cite{Par71} for $k$ a commutative
ring with trivial Picard group.  Later,  
Hopf $H$-Galois extensions \cite{KT,DT} and Hopf subalgebras \cite{OS73,S} 
 have been shown to be  Frobenius extensions of the first and the second 
kinds for $k$ a commutative ring (with a proviso that a Hopf subalgebra
be a $k$-direct summand or pure $k$-submodule in~$H$).  

Although quantum groups, being deformations of the universal enveloping
algebras or the algebra of polynomial functions on Lie groups,
 have  been studied as Hopf algebras
 over fields, we would expect that any
 study of the deformations of affine group
schemes would naturally involve Hopf algebras over commutative rings \cite{WCW}. 

We began a study of a Hopf algebra $H$
over commutative ring $k$ from the point of view
of Frobenius algebras and extensions in \cite{KS} 
starting with previous results in
\cite{LS,Par71,BFS2}. In \cite{KS} we studied a certain class
of Hopf algebras called FH-algebras 
under the condition that $H$ is a Frobenius algebra.  
A purely Frobenius approach to proving the Radford formula
for the fourth power of the antipode $S: H \to H$ 
was taken there.  In this paper we use this approach  to the
Radford formula for a general $H$.   
 The idea of this proof
in \cite{KS} and in the present paper 
 is the following conceptually. First,
from a complete set of  Frobenius data called a Frobenius (coordinate) system for a 
Hopf algebra, we obtain another Frobenius system by applying the 
antipodal anti-automorphism.
Second, we obtain two Nakayama automorphisms with formulas involving $S^{\pm 2}$
 acted on from the right and left, respectively, by  
the left modular function for $H$. Third, the 
 principle that any two  Frobenius systems are unique up to 
an invertible element, called the (Radon-Nikodym) derivative,  
leads after a computation to the modular function for  $H^*$, $b \in H$  as derivative.
Finally, since the two Nakayama automorphisms are related by an inner automorphism 
determined by the derivative, 
we arrive at a conceptually simplified proof for
 the Radford formula for $S^4$. In principle, this technique might produce
nice formulas or new proofs wherever one deals with examples of Frobenius algebras
or extensions. 

In this paper we  will see that a good  working principle 
is that a general Hopf algebra $H$  is very close to being 
an FH-algebra \cite{Par72}.  As noted above,
our main example of this principle 
is to make a Frobenius proof of Radford's formula work for a general finite projective
Hopf algebra $H$. The first part of our paper
is organized around this task as follows. In the Section~2, we present preliminary 
material on a general theory of
$P$-Frobenius algebras \cite{Mor65,Mor67,Par73}
with Frobenius homomorphism, dual bases and 
Nakayama automorphisms, which we also call  {\it a  
Frobenius system} for $H$. To this we add the conceptually useful
{\it comparison theorem} and  {\it transformation
theorem} for $P$-Frobenius algebras. In Section~3 we continue
a review of preliminaries with 
the basic integral theory for a finite projective
Hopf algebra $H$ the conclusion of which is that $H$ is a $P$-Frobenius algebra
with  Frobenius homomorphism $\psi$ very similar to a left integral and
dual bases determined by a left norm $N$. 
In Section~4,  we face the problem that for $H$ the usual definition of modular function
does not work:  the usual definition  
depends on the norm
element being a free generator of the space of integrals, but 
 the space of integrals in $H^*$ is not freely generated by
a left norm.  We instead define a modular function as
the Nakayama automorphism composed with the counit  \cite{KS},
and prove that this plays a successful role. 
In Section~5, we find a formula for the  Nakayama automorphism
of $H$, similar to the formulas in \cite{OS73,FMS}, which eventually leads to 
the proof of Radford's formula in this general case.    Then 
we transform 
the $P$-Frobenius system
for $H$ by the antipode $S$ and prove that the derivative is proportional to 
the distinguished group-like $b \in H$. We finally  
apply the comparison theorem and obtain a complete but conceptual  
proof of Radford's formula for $S^4$ 
in this  general case. 

The rest of the paper is organized as follows. In Section 6,
we show that a finite projective Hopf algebra $H$ is
separable precisely when the counit of its norm is invertible
in some generalized sense for modules.
 We show that if $H$ is separable and 
involutive, then it is strongly separable in Kanzaki's sense;
conversely, as a corollary
of Etingof and Gelaki \cite{EG}, if $H$ is separable and coseparable, it is involutive
(given that $2$ is not a zero-divisor in $k$).   
In Section~7, 
we show that a  Hopf subalgebra pair forms a Frobenius extension
of a {\it third kind}, which is an exotic generalization of Frobenius extensions
of the second kind \cite{NT} and the $P$-Frobenius algebras of Section~2.  
This kind of Frobenius extension 
depends not only on a relative Nakayama automorphism but also on the two Picard
group elements of $k$ represented by the space of integrals of $K^*$ and $H^*$.  
The relative homological algebra of 
Frobenius extensions \cite{Kas2,O} of the first kind 
and Frobenius extensions of the second  and third kinds differs
only in that the functors of co-induction
and induction are naturally equivalent for the first kind
and  differ by a Morita auto-equivalence of the module category 
of the subalgebra 
 for the second and third kinds.
In Section~8, we return to the idea
that a finite projective Hopf subalgebra $H$ is close to being an FH-algebra
by proving that $H$ is a Hopf subalgebra
of an FH-algebra in two ways.  First, we prove that the Drinfel'd
 double $D(H)$ is an FH-algebra. Second,
we find a ring extension $k \subset K$ such that $Pic(K) = 0$:
therefore the FH-algebra $H\otimes_k K$ is a flat extension of $H$.

\section{Preliminaries: $P$-Frobenius Algebras}

In this section we sketch the theory of $P$-Frobenius algebras 
 which generalizes  ordinary Frobenius
algebras and will be needed  in the
later sections (except Proposition (\ref{prop-june18})). 
The material in this section is folkloric and
straightforward applications of for example \cite{Mor65, Mor67, Par73, NEFE}. 
We include short proofs since these have not appeared 
in published form.  
The material after and including Theorem (\ref{th-comp}) is however  
somewhat new.   

Let $k$ be a commutative ring throughout this paper. A tensor
$\otimes$ without subscript  means $ \otimes_k$
as will a homomorphism group $\Hom = \Hom_k$. The
$k$-dual of a $k$-module $V$ is denoted by $V^*$.
If $A$ is a $k$-algebra, its $V$-dual $\Hom(A,V)$ has a standard $A$-bimodule structure
given by $(bfc)(a) := f(cab)$ for every $f \in \Hom(A,V), a,b,c \in A$. 

  Let $P$ be
an invertible  $k$-module throughout, i.e.
$P$ is finite projective of constant rank $1$ \cite{Sil}. The functor
represented by $P \otimes -$ is a Morita auto-equivalence of
the category of $k$-modules, denoted by
${\mathcal M}_k$, and $P$ represents
an isomorphism class in the Picard group $Pic(k)$ of $k$ \cite{Bass,Sil}. 
Let $Q$ be its
inverse as an element of $Pic(k)$, so 
$Q \cong P^*$, and both $P \otimes Q \cong k$ and $Q \otimes P \cong k$
are given by canonical isomorphisms $\phi_1$ and $\phi_2$, respectively,
which we choose so that 
 associativity holds
\begin{equation}
(qp)q' = q(pq')
\label{eq:assoc}
\end{equation}
for every $p \in P$ and $q,q' \in Q$, 
and a corresponding associativity equation on $P\otimes Q \otimes P$ \cite{Bass},
where the values of these isomorphisms are denoted simply by
$p \otimes q \mapsto pq$ and $q \otimes p \mapsto qp$.  Since $\phi_2 \circ \flip \circ
\phi^{-1}_1$
is an automorphism of $k$,
where $\flip: P \otimes Q \to Q \otimes P$
is the ordinary twist map, we have $\chi, \gamma \in k$ such that $\chi \gamma = 1_k$ and
\begin{eqnarray}
pq& = & \gamma qp \nonumber \\
qp & = & \chi pq 
\label{eq:numb}
\end{eqnarray}
for every $p \in P,\, q \in Q$. 

\begin{definition}
A $k$-algebra $A$ is said to be a $P$-Frobenius algebra if
\begin{enumerate}
\item $A$ is finite projective as a $k$-module;
\item $A_A \cong \Hom_k(A,P)_A$.
\end{enumerate}
\label{def-PF}
\end{definition}

If $P \cong P'$, then a $P$-Frobenius algebra is also $P'$-Frobenius.
In particular, if $P \cong k$, then a $P$-Frobenius algebra is  an ordinary Frobenius
algebra. Thus there are no nontrivial $P$-Frobenius algebras over ground rings with trivial Picard
group. The following 
converse statement is {\it false}:  if a $P$-Frobenius algebra is also $P'$-Frobenius, then
$P \cong P'$. This may be somewhat surprising if one recalls that the corresponding statement
is true for
$\beta$-Frobenius extensions \cite{NT}. A counterexample is based on the 
Steinitz isomorphism
theorem $A \oplus B \cong R \oplus AB$ for nonzero ideals $A,B$
 in a Dedekind domain $R$ \cite{JM}: 

\begin{proposition}
 Suppose $R$ is a Dedekind domain and $I$ is a non-principal ideal in $R$ such that
$I \cong I^{-1}$.   Let $A := M_2(R)$.  Then 
\begin{equation}
{}_A\Hom_R(A,I) \cong {}_AA.
\end{equation}
\label{prop-june18}
\end{proposition}
\begin{proof}
Let $F$ denote the field of fraction of $R$, and $e_{ij}$ the matrix units
in $A$. We first note that $\Hom_R(A,I) \cong M_2(I)$, since \[
f \mapsto \left(
\begin{array}{cc}
f(e_{11}) & f(e_{12}) \\
f(e_{21}) & f(e_{22}) 
\end{array}
\right)
\]
is a left $A$-isomorphism if we define the left $A$-module structure on $M_2(I)$
by $X \cdot B := B X^t$ for every $B \in M_2(I), X \in A$. 

By the Steinitz isomorphism theorem, $I \oplus I \cong R \oplus R$ as $R$-modules
determined by a matrix $C \in M_2(F)$ as $(x\ y) \mapsto (x\ y)C^t$.
Then the mapping $X \mapsto (CX)^t$ for every $X \in M_2(I)$ determines
an $R$-isomorphism $\Psi: M_2(I) \rightarrow A$.  But
for every $Y \in A$ we have
\[
\Psi(Y\cdot X) = (CXY^t)^t = Y X^t C^t = Y \Psi(X) \]
whence $\Psi$ is a left $A$-module isomorphism as desired.
\end{proof}

 $A$ is of course a well-known example of a Frobenius algebra over $R$.
That it is also an $I$-Frobenius algebra where $I \not \cong
R$ follows directly from Theorem (\ref{th-pf}) below.
  $R$ is for example realized by the ring of integers of an algebraic number field with two element
ideal class group.

Recall that an algebra $A$ is  QF (quasi-Frobenius) in the sense of  M\"{u}ller \cite{Mul}, 
if $A$ is finite projective as a $k$-module, and $A_A$
is isomorphic to a direct summand of the direct sum of
$n$ copies of $A^*_A$,
for $n \geq 1$.  
It follows straightaway from Definition (\ref{def-PF}) that:

\begin{proposition}
A $P$-Frobenius algebra $A$ is a QF algebra.  
\end{proposition}
\begin{proof}
If $P \oplus N \cong k^n$, then \[
A_A \oplus \Hom_k(A,N) \cong nA^{*} \qed \]
\renewcommand{\qed}{}\end{proof}

Recall that a QF ring $A$ is
artinian and injective as a right or left module over itself \cite{Lam}. If $k$
 is an artinian commutative ring, it has trivial Picard group,
so $A$ in the Proposition is  a QF ring if $k$ is a QF ring \cite{Mul, Kas2}.
 
We shall see below that $P$-Frobenius algebras are much closer to being Frobenius algebras
than QF algebras. 

\begin{theorem}
The following conditions on a $k$-algebra $A$ are equivalent:
\begin{enumerate}
\item  $A$ is a $P$-Frobenius algebra;
\item $A_k$ is finite projective and ${}_AA \cong {}_A\Hom_k(A,P)$;
\item there are $\phi \in \Hom_k(A,P)$, $x_1,\ldots,x_n,y_1,\ldots,y_n \in A$
and $q_1,\ldots,q_n \in Q$ such that 
\begin{equation}
\sum_i \phi(ax_i)q_i y_i = a
\label{eq:left}
\end{equation}
or 
\begin{equation}
\sum_i x_i q_i \phi(y_ia) = a
\label{eq:right}
\end{equation}
for every $a \in A$. ($\phi$ is referred to as a Frobenius homomorphism
and $\{ x_i \},\{ q_i\},\{ y_i \}$ as  dual bases for $\phi$.) 
\end{enumerate}
\label{th-pf}
\end{theorem}
\begin{proof}
(1 $\Rightarrow$ 2.)  We compute using the Hom-tensor relation: 
\begin{eqnarray*}
{}_A\Hom_k(A,P) & \cong & \Hom_k(\Hom_k(A,P)_A, P) \\
& \cong & {}_A\Hom_k(A^* \otimes P,P) \\
& \cong & {}_A\Hom_k(A^*,k) \cong {}_AA,
\end{eqnarray*} 
since $P$ is an invertible module.  

(2 $\Rightarrow$ 3.)  Given $\Psi: {}_AA \stackrel{\cong}{\rightarrow}
{}_A\Hom_k(A,P)$
and $\phi := \Psi(1_A)$, then $\Psi(a) = a\phi$ for every $a \in A$.  
Then ${}_AA \otimes Q \cong {}_AA^*$ via $a \otimes q \mapsto a\phi q$.
If $\{ y_i \in A \},\, \{ f_i \in A^* \}$ is a finite projective base for $A_k$,
one finds $x_{ij} \in A$, $q_{ij} \in Q$ such that 
 $\sum_{j} x_{ij} \phi q_{ij} = f_i$. Setting $y_{ij} := y_i$ for each $i$ and $j$, 
we have for every $a \in A$, 
\begin{eqnarray*}
a & = & \sum_i f_i(a) y_i \\
& = & \sum_{i,j} (x_{ij} \phi)(a)q_{ij} y_{ij}  = \sum_{i,j} \phi(ax_{ij})q_{ij} y_{ij}.
\end{eqnarray*}
We merely reindex to get Eq.\ (\ref{eq:left}). 
Eq.\ (\ref{eq:right}) follows from a computation showing
 $\Psi(\sum_i x_i q_i \phi(y_ia))(x) = \Psi(a)(x)$ for $x,a \in A$,
which is similar to \cite[1.3]{FMS}.

(3 $\Rightarrow$ 1.)  Suppose $\sum_{i=1}^n x_i q_i(\phi  y_i ) = \id_A$.  Then 
$A$ is finite projective.
Define $\Psi: A_A \rightarrow \Hom_k(A,P)_A$ by $\Psi(a):= \phi a$ for every
$a \in A$.  Then $\Psi $ is epi since for every $f \in \Hom_k(A,P)$
we have $\Psi(\sum_i f(x_i)q_i y_i)(a) = f(a)$ for every $a \in A$. 
Since $\Psi: A \rightarrow \Hom_k(A,P) \cong A^* \otimes P$ 
is an epimorphism between finite projective modules of the same local rank,
(i.e. ${\mathcal P}$-rank
for every prime ideal ${\mathcal P}$ in $k$), $\Psi$ is bijective
\cite{Sil,Par71}.  

A similar argument shows that we may establish 
Condition 2 from Eq.\ (\ref{eq:left}).
\end{proof}

Throughout this section, we continue our use of the notation $\phi$ and $x_i,q_i,y_i$ for the Frobenius
homomophism and dual base of a $P$-Frobenius algebra $A$. A QF ring has a Nakayama
permutation on the set of simples modules induced by taking the socle of 
the corresponding projective indecomposable modules 
\cite{Lam}.  Frobenius algebras moreover have Nakayama {\it automorphisms} \cite{Kas2}.  
We  next see that $P$-Frobenius algebras also have Nakayama
automorphisms.  
\begin{corollary}
In a $P$-Frobenius algebra
$A$ there is an algebra automorphism $\nu: A \rightarrow A$ given by
\begin{equation}
a\phi = \phi \nu(a)
\label{eq:naka}
\end{equation}
for every $a \in A$.  (Call $\nu$ the Nakayama automorphism.)
\label{cor-Naka}
\end{corollary}
\begin{proof}
In the proof of the last theorem we established 3 $\Rightarrow$ 1 by
showing $a \mapsto \phi a$, for every $a \in A$, is an isomorphism. 
As we noted, we may equally well establish 3 $\Rightarrow$ 2 in this proof 
by showing that $a \mapsto a \phi$ is an isomorphism ${}_A A \cong {}_A\Hom_k(A,P)$.
Since $a\phi \in \Hom_k(A,P)$ for each $a \in A$, it follows that there is a unique $a' \in A$
such that $a\phi = \phi a'$.  One defines $\nu(a) = a'$ and easily checks that
$\nu$ is an automorphism. 
\end{proof}

In this respect a $P$-Frobenius algebra is almost Frobenius:
of course, $\nu$ measures the deviation
of $\phi$ from satisfying the trace condition $\phi(ab) = \phi(ba)$
for every $a,b \in A$. If $\nu$ is inner, $A$ will possess such a trace-like
Frobenius homomorphism and is called a {\it symmetric
$P$-Frobenius algebra}. 
We fix the data $(\phi,x_i,q_i,y_i,\nu)$ for the rest of this section
and refer to this as the {\it Frobenius system} of $A$
 in this paper.

\begin{proposition}
Given a $P$-Frobenius algebra $A$, the dual base tensor $\sum_i x_i \otimes q_i \otimes y_i$
satisfies $\forall a \in A$:
\begin{enumerate}
\item $\sum_i ax_i \otimes q_i \otimes y_i = \sum_i x_i \otimes q_i \otimes y_ia$, and 
\item $\sum_i x_ia \otimes q_i \otimes y_i = \sum_i x_i \otimes q_i \otimes \nu(a) y_i$.
\end{enumerate}
\label{prop-enumb}
\end{proposition}
\begin{proof}
We give only the proof of the second equation, the first being similar.
By Eqs.\ (\ref{eq:right}), (\ref{eq:assoc}), (\ref{eq:naka}) and (\ref{eq:left}),
we compute:
\begin{eqnarray*}
 \sum_i x_ia \otimes q_i \otimes y_i & = & \sum_{i,j} x_jq_j\phi(y_j x_i a) \otimes q_i \otimes y_i \\
& = & \sum_{i,j} x_j \otimes q_j \otimes \phi(y_jx_i a)q_i y_i \\
& = & \sum_{i,j} x_j \otimes q_j \otimes \phi(\nu(a) y_jx_i)q_iy_i \\
& = & \sum_j x_j \otimes q_j \otimes \nu(a) y_j \qed
\end{eqnarray*}
\renewcommand{\qed}{}\end{proof}

We next  prove  that 
$P$-Frobenius systems for $A$ are unique up to an invertible element in $A$,
which we call the {\it comparison theorem}.  

\begin{theorem}(``Comparison Theorem'').
Suppose $(\phi,x_i,q_i,y_i)$
and $(\phi',x_j',q_j',y_j')$ are two $P$-Frobenius systems for a 
$P$-Frobenius algebra $A$.  Then there is $d \in A^{\circ}$
such that 
\begin{equation}
\phi' = \phi d
\label{eq:dee}
\end{equation}
 and 
\begin{equation}
\sum_j x_j' \otimes q_j' \otimes y_j' = \sum_i x_i \otimes q_i \otimes d^{-1}y_i.
\label{eq:deeday}
\end{equation}
If $\nu, \nu'$ are the Nakayama automorphisms of $\phi$ and $\phi'$, then
$\forall a \in A$, 
\begin{equation}
v'(a) = d^{-1}v(a) d.
\label{eq:aprilsun}
\end{equation}
\label{th-comp}
\end{theorem}
\begin{proof}
Since $\phi$ and $\phi'$ freely generate $\Hom_k(A,P)$ as right $A$-modules,
Eq.\ (\ref{eq:dee}) is clear with $d$ an invertible in $A$.

To verify Eq.\ (\ref{eq:deeday}), we note that
\begin{equation}
\sum_i x_i q_i \phi d (d^{-1}y_i a) = a 
\label{eq:inter}
\end{equation}
for every $a \in A$.  There is an isomorphism \[
A \otimes Q \otimes
 A \cong \End_k(A)
\]
given by $a \otimes q \otimes b \mapsto aq\phi'b$, for every $a,b \in A,
q \in Q$, since $A \otimes A^* \cong \End_k(A)$
and $Q \otimes A \cong A^* $. 
 Eq.\ (\ref{eq:deeday}) follows from the injectivity of this
mapping and Eq.\ (\ref{eq:inter}).

We note that for every $x,a \in A$
\begin{equation}
\phi'(xa) = \phi'(\nu'(a)x) \ \Leftrightarrow \  \phi(dxa) = \phi(\nu(a)dx)
= \phi(d\nu'(a)x)
\end{equation}
the last equation implying that for all $a \in A$, 
\[
\nu(a) d = d \nu'(a)
\]
which is equivalent to Eq.\ (\ref{eq:aprilsun}). 
\end{proof}

We also need to know the effect of an algebra anti-automorphism on a Frobenius
system, as given in the following {\it transformation theorem.}

\begin{theorem}(``Transformation Theorem'').
Let $A$ be a $P$-Frobenius algebra with Frobenius system $(\phi,x_i,q_i,y_i,\nu)$.
If $\alpha$ is a $k$-algebra anti-automorphism of $A$, then 
\begin{equation}
(\alpha \phi,\ \chi \overline{\alpha}(y_i),\ q_i,\ \overline{\alpha}(x_i), \
\overline{\alpha}\circ \overline{\nu} \circ \alpha)
\end{equation}
is another Frobenius system for $A$, 
where $\overline{\alpha}$ and $\overline{\nu}$
denote the inverses of $\alpha$ and $\nu$, and $\alpha \phi := \phi \circ \alpha$.
\label{prop-antitransform}
\end{theorem}
\begin{proof}
We compute using the identity $\alpha(ab) = \alpha(b) \alpha(a)$ for all $a,b \in A$:
\[
a = \sum_i x_i q_i \phi(y_i a) =
 \sum_i \chi (\alpha \phi)(\overline{\alpha}(a)\overline{\alpha}(y_i)) q_i x_i,
\] 
and by applying $\overline{\alpha}$ to both sides we obtain
\[
\overline{\alpha}(a) =  \sum_i \chi (\alpha \phi)(\overline{\alpha}(a)\overline{\alpha}(y_i)) q_i 
\overline{\alpha}(x_i).
\]
It follows from Theorem (\ref{th-pf}) that $\alpha \phi$ is a Frobenius homomorphism
with dual bases $\{  \chi \overline{\alpha}(y_i) \}$, $\{ q_i\} $, $ \{ 
 \overline{\alpha}(x_i) \}$.  

We compute the Nakayama automorphism $\eta$ for $\alpha \phi$ in terms of $\alpha$ and
$\nu$: for all $a,b \in A$, 
\[
\phi(\alpha(a)\alpha(b)) = (\alpha \phi)(ba) = (\alpha \phi)(\eta(a)b) 
= \phi(\alpha(b) \alpha \eta(a))
= \phi( (\nu \alpha \eta)(a) \alpha(b))
\]
by applying Eq.\ (\ref{eq:naka}) twice.  Since $\phi$ freely generates $A^*$,
it follows that $\nu \circ \alpha \circ \eta = \alpha$, whence 
\begin{equation}
\eta = \overline{\alpha} \circ \overline{\nu} \circ \alpha. \qed
\label{eq:eta}
\end{equation}
\renewcommand{\qed}{}\end{proof}

We will need the following lemma in our last section.  
\begin{lemma}
If $A$ is a $P$-Frobenius algebra and $B$ is a $Q$-Frobenius algebra,
then the tensor product algebra $A \otimes B$ is a $P \otimes Q$-Frobenius
algebra.
\label{lemma-june1}
\end{lemma}
\begin{proof}
First, $ C := A \otimes B$ is finite projective as a $k$-module.
Secondly,
\[
{}_CC \cong {}_AA \otimes {}_BB \cong {}_A\Hom(A,P) \otimes {}_B\Hom(B,Q) \cong {}_C\Hom(C,P \otimes Q),
\]
since $A$, $B$, $P$ and $Q$ are finite projective $k$-modules.
\end{proof}

\section{Preliminaries II:
 Hopf Algebras as P-Frobenius Algebras}

Let $H$ be a Hopf algebra over a commutative ring $k$, 
which is finite (i.e., finitely generated) projective
as a $k$-module, throughout this paper unless otherwise stated.
In this section, we review  
 the Hopf module structure on the dual Hopf algebra $H^*$ \cite{LS, Par71}
and the $P$-Frobenius structure on $H$ \cite{Par73}.  
For the convenience of the reader
we   offer  proofs for the propositions that have not been published.

For the Hopf algebra $H$ we denote its comultiplication by $\Delta: H \rightarrow H \otimes H$,
its counit by $\epsilon$, and its antipode by $S$.  The values of $\Delta$ are denoted by
$\Delta(x) = \sum x\1 \otimes x\2$.  If $M$ is a right comodule over $H$ the values of
its coaction on an element $m \in M$ is  denoted by $\sum m\0 \otimes m\1$. The dual
of $H$ is itself a Hopf algebra $H^*$ where its multiplication is the convolution product (dual
to $\Delta$),
comultiplication is the dual of multiplication
on $H$, the counit is $1 \in H \cong H^{**}$ ($x \mapsto $ evaluation at $x$).
We also denote its antipode by $S$ where the context is clear.

\begin{proposition}
If $H$ is a finite projective Hopf algebra, then $H^*$ is right Hopf module.
\label{prop-one}
\end{proposition}
\begin{proof}[Sketch of Proof in \cite{LS,Par71}]
The natural left $H^*$-module structure
on the dual algebra $H^*$ induces a comodule structure mapping $\chi: H^* \rightarrow H^* \otimes H$,
determined by
\begin{equation}
gh = \sum h\0 g(h\1)
\label{eq:ordinaire}
\end{equation}
for every $g,h \in H^*$. 
The right $H$-module structure on $H^*$ is given by $(h^* \cdot h)(x) = h^*(xS(h))$
for every $x,h \in H$ and $h^* \in H^*$.  A rather long computation 
shows this compatible with the $H^*$-comodule structure in the sense of Hopf modules.
\end{proof}

\begin{proposition}
A right Hopf module $M$ over a finite projective Hopf algebra $H$
 is isomorphic to the trivial Hopf module,
$M \cong P(M) \otimes H$, where 
\[
P(M) = \{ m \in M | \, \chi(m) = m \otimes 1_H \}
\]
is a $k$-direct summand of $M$ and $\chi: M \rightarrow M \otimes H$ denotes
the right $H$-comodule structure mapping.
\label{prop-two}
\end{proposition}
\begin{proof}[Sketch of Proof in \cite{Par71}]
One shows that the map $M \rightarrow M$
given by $m \mapsto \sum S(m\0) m\1$ is a $k$-linear projection onto $P(M)$.
Then the mapping $\beta: M \rightarrow P(M) \otimes H$ given
by $\beta(m) = \sum m\0 S(m\1) \otimes m\2$
has inverse given by the Hopf module map $\alpha: P(M) \otimes H \rightarrow M$
defined by $\alpha(m \otimes h) = mh$. 
\end{proof} 

\begin{corollary}
The $k$-module $P(H^*)$ associated to a Hopf algebra $H$
by Propositions (\ref{prop-one}) and (\ref{prop-two})
 is an invertible  $k$-direct summand in $H^*$.
\label{cor-sil}
\end{corollary}
\begin{proof}
Since $P(H^*) \otimes H \cong H^*$ and $H$, $H^*$ have the same local ranks, 
 it follows that the finite projective
$k$-module 
$P(H^*)$ has constant rank $1$.  Then  $P(H^*) \otimes P(H^*)^* \cong k$ and 
$P(H^*)$ is invertible 
 \cite{Sil}.
\end{proof}

We  note that $P(H^*)$ is the space of left integrals $\int^{\ell}_{H^*}$
in $H^*$:
\begin{equation}
P(H^*) = \{ f \in H^* | gf = g(1)f \}
\end{equation}
which follows from Eq.\ (\ref{eq:ordinaire}) since $\sum f\0 \otimes f\1 = f \otimes 1$.
 
\begin{proposition}
The antipode $S$ of a finite projective Hopf algebra $H$ is bijective.
\end{proposition}
\begin{proof}[Sketch of Proof in \cite{Par71}]
Assuming that $S(x) = 0$, one
then notes that multiplication from the right by $x$ on $P(H^*) \otimes H$ is zero
by the existence of the ($H$-module) isomorphism $\alpha: P(H^*) \otimes H \rightarrow H^*$
in Proposition (\ref{prop-two}).
If $k$ is field $P(H^*) \cong k$ and it is clear that $x $ is then zero.
The general case follows from a localization argument.  Surjectivity for $S$
is apparent if $k$ is a field, and the general case follows again from a localization
argument. 
\end{proof}

Denote the composition-inverse of $S$ by $\es$. 

\begin{proposition}[\cite{Par73}] 
If $H$ is a finite projective Hopf algebra and \[
P := P(H^*)^*,
\]
 then $H$ is a $P$-Frobenius algebra.
\label{prop-pa}
\end{proposition}
\begin{proof}
We set $\Phi: P(H^*) \otimes H \stackrel{\cong}{\longrightarrow} H^*$, $f \otimes x \mapsto
f \cdot x$, where we note that the right $H$-module structure is related to
the standard left $H$-module structure on $H^*$ via a twist by $S$: for every $g \in H^*,
x,y \in $ \[
(g\cdot x)(y) = g(y S(x)) = (S(x)g)(y).\]
  Let \(
 Q : = P(H^*),
\)
which is canonically isomorphic to the dual of $P$, and satisfies $P \otimes Q \cong k$
by Corollary (\ref{cor-sil}).

Define $\Psi': H \rightarrow \Hom_k(H,P)$ as the composite of the right $H$-module
isomorphisms
\[
H \longrightarrow P \otimes Q \otimes H \stackrel{1 \otimes \Phi}{\longrightarrow} 
P \otimes H^* \longrightarrow \Hom_k(H,P).
\]
It is easy to check that 
\begin{equation}
\Psi'(x)(y)(q) := \Phi(q \otimes x)(y) = q(yS(x))
\label{eq:icebox}
\end{equation}
for all $x,y \in H$ and $q \in Q$. 

Now let $\Psi := \Psi'\circ \es$.  $\Psi$ is a Frobenius isomorphism ${}_H H \cong {}_H \Hom_k(H,P)$,
since $\es$ is an anti-automorphism of $H$ and 
\[
\Psi(xy) = \Psi'(\es(y)\es(x)) = \Psi(y)\cdot \es(x) = x\Psi(y). \qed
\]
\renewcommand{\qed}{}\end{proof}

Gabriel has an example of a finite projective Hopf algebra which is {\it not} a Frobenius
algebra \cite{Par71}. 

\begin{corollary}
The Frobenius homomorphism $\psi: H \rightarrow P$ defined by
the theorem satisfies for every $a \in H$ 
\begin{equation}
\sum a\1 \otimes \psi(a\2) = 1 \otimes \psi(a)
\label{eq:li}
\end{equation}
\end{corollary}
\begin{proof}
We note that the Frobenius homomorphism
$\psi:=\Psi(1) = \Psi'(1)$ satisfies by Eq.\ (\ref{eq:icebox}),
for  every $q \in P(H^*), a \in H$, \[
\psi(a)(q) = q(a),
\] and \[
q(a)1_H = \sum a\1 q(a\2)
\]
since $q \in \int^{\ell}_{H^*}$. 

Since $P = Q^*$ and $H$ is finite projective over $k$,
we canonically identify $H \otimes P \cong \Hom_k(Q,H)$,
and compute $\forall q \in Q, a \in H$
\[
(\sum a\1 \otimes \psi(a\2))(q) = \sum a\1 \psi(a\2)(q) = \sum a\1 q(a\2) = 1_H q(a) = (1 \otimes \psi(a))(q)
\]
whence Eq.\ (\ref{eq:li}). 
\end{proof}

If $\int^{\ell}_{H^*} \cong k$, we see from the theorem and the corollary that $H$
is an ordinary Frobenius algebra with Frobenius homomorphism a left integral in $H^*$:
this is called an {\it FH-algebra}  \cite{Par72, KS}. 
Conversely, we have
the following result. 

\begin{proposition}
If $H$ is a Frobenius algebra and Hopf algebra, then $H$ is an FH-algebra.
\end{proposition}
\begin{proof}
We use the fact that the $k$-submodule of integrals of an augmented Frobenius
algebra is free of rank $1$ (cf.\ Lemma (\ref{lemma-D}), 
\cite[Theorem 3]{Par71} or \cite[Prop.\ 3.1]{KS}).  Then
$\int^{\ell}_H \cong k$.  It follows from Proposition (\ref{prop-pa})
 that the dual
 Hopf algebra $H^*$ is a Frobenius algebra. Whence $\int^{\ell}_{H^*} \cong k$
and $H$ is an FH-algebra. 
\end{proof}
 
Next 
we  obtain  as in \cite{Par73} a left norm for the Frobenius homomorphism $\psi: H \rightarrow P$
and study its properties.  
Since $x \mapsto x\psi$ is an isomorphism ${}_H H \rightarrow {}_H\Hom(H,P)$
and $\Hom(H,P) \otimes Q \cong H^*$ affords a canonical identification,
it follows that there are elements $N_i \in H, q_i \in Q$ such that the counit of $H$, 
\begin{equation}
\epsilon \stackrel{\cong}{\longmapsto} \sum_i N_i \psi \otimes q_i .
\label{eq:ln}
\end{equation}
Call $N := \sum_i N_i \otimes q_i$ in $H \otimes Q$ the {\it left norm} of $\psi$,
and note that $\sum_i \psi(aN_i)q_i = \epsilon(a)$ for every $a\in H$. 
In the natural left $H$-module ${}_H H \otimes Q$ we have
\begin{equation}
aN = \epsilon(a)N,
\label{eq:lemA}
\end{equation}
since both $aN$ and $\epsilon(a)N$ map to $\epsilon(a)\epsilon$ under
the composite isomorphism, $H \otimes Q \stackrel{\cong}{\rightarrow} \Hom_k(H,P) \otimes Q
\stackrel{\cong}{\rightarrow} H^*$ given by $a \otimes q \mapsto a \psi q$.

For all $p \in P$, we note that
\begin{equation}
\sum_i N_i q_i(p) \in \int^{\ell}_H,
\label{eq:lemB}
\end{equation}
since this follows by applying Eq.\ (\ref{eq:lemA}) to $p$. 

\begin{proposition}[\cite{Par73}] 
If $H$ is a Hopf algebra with Frobenius homomorphism $\psi$ given
above and left norm $\sum_i N_i \otimes q_i$, then
the dual bases  for $\psi$ is given by 
\begin{equation}
\{ N\II \},\, \{ q_i\} ,\, \{ \es(N\I) \}
\end{equation}
\label{prop-par73}
\end{proposition}
\begin{proof}
We compute
as in \cite[Lemma 3.16]{Par73},  using Eq.\ (\ref{eq:li}) at first
and Eq.\ (\ref{eq:lemA}) next 
 (for every $ a \in A$):
\begin{eqnarray*}
\sum \psi(aN\II)q_i \es(N\I) & = & \sum a\1N\II (\psi(a\2N\III)q_i) \es(N\I) \\
& = & \sum a\1 \psi(a\2N_i)q_i \\
& = & \sum a\1 \epsilon(a\2) \psi(N_i)q_i = a\epsilon(1) = a.
\end{eqnarray*}
It follows from Theorem (\ref{th-pf}) that $\{ N\II \},\, \{ q_i\} ,\, \{ \es(N\I) \}$
are dual bases for $\psi$.
\end{proof}

\section{Pinning Down the Modular Functions}

In this section
we give a definition of modular function in Eq.\ (\ref{eq:em}) based on \cite{KS},
and find two formulas, Eqs.\ (\ref{eq:lemC}) and (\ref{eq:coldsun}) which will be used later.
The rest of this section is somewhat technical and might be browsed on a first reading. 

It follows from appying $S$ to the equation in the last proof, and setting $a = 1$, that
\begin{equation}
\sum_i (\psi q_i) \rightharpoonup N_i = 1,
\label{eq:seagulls}
\end{equation}
where $\psi q_i \in H^*$ is the mapping $a \mapsto \psi(a)q_i$ for each $i$ and $a\in H$.
Of course $1 \in H^{**} \cong H$ is the counit of $H^*$. 
It follows from Eqs.\ (\ref{eq:ln}) and (\ref{eq:seagulls})
that the antipode on the dual Hopf algebra $H^*$ is given
by
\begin{equation}
S(g) = \sum N_i(g (\psi q_i)\2) (\psi q_i)\1,
\label{eq:donnerstag}
\end{equation}
since one computes that $\sum g\1 S(g\2) = g(1)\epsilon$ for every $g \in H^*$.  

\begin{proposition}
 $H$ is a Hopf algebra and $P$-Frobenius algebra if and only if 
$H^*$ is a Hopf algebra and $P^*$-Frobenius algebra.
\label{prop-fixed}
\end{proposition}
\begin{proof}
Let $Q = P^*$. 
It suffices to show the forward implication.  Let $p_i \in P$ be 
such that $\sum_i q_i p_i = 1_k$.  Then  Eq.\ (\ref{eq:donnerstag})
implies that 
\begin{equation}
(N_i \otimes q_i,(\psi q_i)\2,\  p_i, \ \es(\psi q_i)\1)
\end{equation}
is a $Q$-Frobenius system for $H^*$, where we identify $H \otimes Q \cong \Hom(H^*,Q)$
via the obvious isomorphism. 
\end{proof}

We next define a {\it left modular function}
for a Hopf algebra $H$.  We continue the notation established in the previous section. 

\begin{definition}
Define the {\it left modular function}, or {\it 
left distinguished group-like element}, $m: H \rightarrow k$ by
\begin{equation}
m := \epsilon \circ \nu
\label{eq:em}
\end{equation}
where $\nu$ is the Nakayama automorphism of $H$ relative
to $\psi$ (cf.\ Corollary (\ref{cor-Naka})).
\end{definition}

First note that $m$ does not depend on the choice of Nakayama automorphism,
since $\epsilon(d \nu(a) d^{-1}) = \epsilon(\nu(a))$ for every $a \in A$. 
Next note that $m$ is an algebra homomorphism (an augmentation in fact),
and therefore a group-like element in the dual Hopf algebra $H^*$. 
With respect to the natural right $H$-module
$H_H \otimes_k Q$, we note that for all $a \in H$,
\begin{equation}
Na = Nm(a),
\label{eq:lemC}
\end{equation}
since $Na$ is mapped into
$\sum_i N_i a \psi \otimes q_i = \sum_i N_i \psi \nu(a) \otimes q_i$, 
then into $\epsilon(\nu(a)) \epsilon = m(a) \epsilon$,
under the canonical isomorphism $H \otimes Q \cong  H^*$. 

Let $A$ be an algebra with augmentation $\epsilon$,
 ${}_AM_A$ an $A$-bimodule and define
the $k$-module of left integrals in $M$ as $\int^{\ell}_M := \{ x \in M| ax = \epsilon(a)x \}$.
For a Hopf algebra and $P$-Frobenius algebra $H$ we consider the natural
$H$-bimodule ${}_HH_H \otimes Q$ in the lemma below.

\begin{lemma}
Given Hopf algebra $H$ and Frobenius homomorphism $\psi$, $\int^{\ell}_{H \otimes Q}$
is a sub-bimodule freely generated by 
the left norm $N = \sum_i N_i \otimes q_i$ and a $k$-direct summand of $H \otimes Q$.
\label{lemma-D}
\end{lemma}
\begin{proof}
$N$ is left integral by Eq.\ (\ref{eq:lemA}).  
We recall the isomorphism $H \otimes Q \stackrel{\cong}{\rightarrow} \Hom(H,P) \otimes Q
\stackrel{\cong}{\rightarrow} H^*$ given by $a \otimes q \mapsto
(a\psi)q$. Given
$T = \sum_i T_i \otimes q'_i \in \int^{\ell}_{H \otimes Q}$, denote $\phi(T) := \sum_i \psi(T_i)q'_i
\in k$, and note that, for all $x \in H$,
\[
\sum_i \psi(xT_i)q'_i = \epsilon(x) \phi(T) = \sum_i \psi(xN_i)q_i \phi(T).
\]
Whence 
\begin{equation}
T = \phi(T)N.
\label{eq:useful}
\end{equation}
Thus, $N$ generates $\int^{\ell}_{H \otimes Q}$ and the mapping of $H \otimes Q \rightarrow 
\int^{\ell}_{H \otimes Q}$ given by $x \otimes q \mapsto \psi(x)q N$ is a $k$-linear
projection.

If $\lambda \in k$ such that $\lambda N = 0$, then
\[
0 = \sum_i \psi(N_i)q_i \lambda = \epsilon(1) \lambda = \lambda,
\]
so $N$ freely generates $\int^{\ell}_{H \otimes Q}$. 
\end{proof}

We similarly define right integrals in a bimodule over an augmented algebra,
and prove a right-handed version of the lemma. It follows from
Lemma (\ref{lemma-D}) that $T := \sum_i \es(N_i) \otimes q_i$ is a right integral
that freely generates $\int^r_{H \otimes Q}$,
since $\epsilon \circ \es = \epsilon$ and $\es$ is an anti-automorphism of $H$. 
By Proposition (\ref{prop-par73}), we compute
\begin{eqnarray*}
T = \sum_{i,j,(N_i)} \psi(\es(N_j)N\II)q_i \es(N\I) \otimes q_j & =
& \sum \psi(\es(N_j))q_i \es(N\I)\epsilon(N\II) \otimes q_j \\
&= & T (\sum_j q_j \psi(\es(N_j))),
\end{eqnarray*}
whence
\begin{equation}
\sum_j q_j \psi(\es(N_j))) = 1_k.
\label{eq:coldsun}
\end{equation}
It follows that $T$ is a right norm in the sense that $\sum_i q_i \psi \es(N_i) = \epsilon$.

\begin{lemma}
$\psi$ and $\psi \circ \es$ are left and right norms in the natural $H$-bimodule
$H^* \otimes P \cong \Hom_k(H,P)$.
\label{lemma-E}
\end{lemma}
\begin{proof}
Proposition (\ref{prop-fixed})
shows that $N \in H \otimes Q$ is a Frobenius homomorphism for the dual Hopf algebra $H^*$.
The concepts of left and right norm relative to $N$ make
sense in the $H^*$-bimodule $H^* \otimes Q$. But Eq.\ (\ref{eq:seagulls}) implies that 
 $\psi \in \Hom(H,P) \cong H^* \otimes P$ is a left norm for $N$.
Similarly, $\sum_i S(N_i) \otimes q_i$ is a Frobenius homomorphism $H^* \rightarrow Q$ 
by applying the anti-automorphism $S$ as in Theorem (\ref{prop-antitransform}), and 
 $\es \psi$ is a right norm. 
\end{proof}

One easily checks that $\sum_i S(N_i) \otimes q_i$ is a right norm in $H \otimes Q$
for $\psi \circ \es$. Since $H^*$ is a $Q$-Frobenius algebra, it has a Nakayama
automorphism $\nu^*$, which we make formal use of  below.  

\begin{definition}
Let $b \in  H$, where $H$ is canonically identified with $H^{**}$, be the left modular function
defined by
\begin{equation}
b = \eta \circ \nu^*
\end{equation}
where $\eta$ is the counit of $H^*$ defined by $\eta(f) =  f(1)$ for every $f \in H^*$.
\end{definition}

It follows from Eq.\ (\ref{eq:lemC}) and Lemma (\ref{lemma-E}) that for every
$f \in H^*$, 
\begin{equation}
\psi f = \psi f(b),
\label{eq:defB}
\end{equation}
where $\psi \in H^* \otimes P$ has the natural $H^*$-bimodule structure. 

\section{An Application to Radford's Formula} 

We now compute a formula for the Nakayama automorphism of $\psi: H \rightarrow P$
in terms of the square
of the antipode and $m$.  The notation $g \rightharpoonup a : = \sum a\1 g(a\2)$
and $a \leftharpoonup g := \sum g(a\1) a\2$ denotes the usual left and right
module actions of the convolution algebra $H^*$ on $H \cong H^{**}$.

\begin{theorem}
The Nakayama automorphism $\nu$ for $\psi: H \rightarrow P$ is given by
\begin{equation}
\nu(a) = \es^2(m \rightharpoonup a) = m \rightharpoonup \es^2(a)
\label{eq:nu}
\end{equation}
\label{th-C}
\end{theorem}
\begin{proof}
The rightmost equation follows from noting that $m$ is a group-like element in $H^*$,
whence $m \circ S = m^{-1}$ and $m \circ S^2 = m$: i.e., $S^2$ and $\es^2$ fix $m$.

The leftmost equation is computed 
below and follows \cite[Satz 3.17]{Par73} until (\ref{eq:nolonger}): for every $a\in H$, 
\begin{eqnarray}
S^2(\nu(a)) &  = & S^2(\sum\psi(N\II a)q_i \es(N\I)) \nonumber \\
& = & \sum S(N\I) \psi(N\II a) q_i \nonumber \\
& = & \sum S(N\I)N\II a\1 \psi(N\III a\2)q_i \nonumber \\
& = & \sum a\1 \psi(N_i a\2) q_i 
\label{eq:nolonger} \\
& = & \sum a\1 m(a\2) \psi(N_i) q_i \nonumber \\
& = & m \rightharpoonup a \nonumber
\end{eqnarray}
by Eqs.\ (\ref{eq:naka}), (\ref{eq:li}), (\ref{eq:lemC}) and (\ref{eq:ln}),
respectively. 
\end{proof}

Since $H$ has Frobenius system $( \psi,\, N\II,\, q_i,\, \es(N\I),\,
\nu)$, it follows from Theorem (\ref{prop-antitransform}) 
that we obtain another Frobenius system
by applying the algebra (and coalgebra) anti-automorphism $\es$:

\begin{proposition}
A Hopf algebra $H$ with left norm $N$  has Frobenius system
\begin{equation}
( \es \psi,\ \chi N\I,\ q_i,\ S(N\II),\ \alpha)
\label{eq:fs}
\end{equation}
where $\es \psi$ satisfies a ``right integral-like equation,''
\begin{equation}
(\es \psi)(x) \otimes 1_H = \sum (\es \psi)(x\1) \otimes x\2
\label{eq:ri}
\end{equation}
and the Nakayama automorphism,
\begin{equation}
\alpha(x) = S^2(x) \leftharpoonup m
\label{eq:alpha}
\end{equation}
for every $x \in H$.
\label{prop-F}
\end{proposition}
\begin{proof}
The dual bases (\ref{eq:fs}) follows directly from Theorems (\ref{prop-antitransform})
and (\ref{prop-par73}). Eq.\ (\ref{eq:ri}) follows from Eq.\ (\ref{eq:li}) since
$\es$ is a coalgebra anti-automorphism.

To compute the Nakayama automorphism we first need to find the inverse
of Eq.\ (\ref{eq:nu}): for all $a\in H$, 
\begin{equation}
\overline{\nu}(a) = S^2(m^{-1} \rightharpoonup a) = m^{-1} \rightharpoonup S^2(a).
\label{eq:nuinv}
\end{equation}
Next we apply Eq.\ (\ref{eq:eta}) where $\es$ is the anti-automorphism: 
\begin{eqnarray*}
\alpha(x) & = &  (S \circ \overline{\nu} \circ \es)(x) \\
& = &  S(m^{-1} \rightharpoonup S(x)) \\
& = &  S(\sum S(x\2) m^{-1}(S(x\1))) \\
& = & S^2(x) \leftharpoonup m ,
\end{eqnarray*}
since $m \circ S = m^{-1}$ and
$S^2$ is an algebra and coalgebra {\it automorphism}.
\end{proof}

By the comparison theorem, we know that the two Frobenius homomorphisms 
$\psi$ and $\es \psi$ are related by an invertible element $d$ called the derivative:
$\es \psi = \psi d$.  The next proposition
shows that $d$ is proportional to the left distinguished group-like element $b$ of $H^*$.
\begin{proposition}
If $\psi$ is a Frobenius homomorphism for the Hopf algebra $H$, then 
\begin{equation}
\psi \circ \es = \gamma \psi b
\end{equation}
\label{prop-G}
\end{proposition}
\begin{proof}
We first show that $\psi b$ is a right integral in the $H^*$-bimodule
$H^* \otimes P$. Recall that $H^* \otimes P $ is canonically identified with $\Hom_k(H,P)$
Let $f \in H^*$, then
\[
(\psi b)f  = [\psi(fb^{-1})]b = [\psi((fb^{-1})(b))]b = (\psi b)f(1)
\]
since $\Delta(b) = b \otimes b$. 

Since $\psi \circ \es$ is a right norm it follows that there is $\lambda \in k$
such that $\psi \circ \es = \lambda (\psi b)$. But
comparing Eq.\ (\ref{eq:coldsun}) to the application below of Eq.\ (\ref{eq:lemA}):
\[
\sum_i q_i (\psi b)(N_i) = \chi \epsilon(b) \epsilon(1) = \chi,
\]
shows that $\lambda = \gamma$ (cf.\ Eq.\ (\ref{eq:numb})). 
\end{proof}

\begin{theorem}
If $H$ is a finite projective Hopf algebra with left distinguished group-like elements
$b \in H$ and $m \in H^*$, then for every $a \in H$,
\begin{equation}
S^4(a) = b^{-1}(m \rightharpoonup a \leftharpoonup m^{-1})b.
\label{eq:radford}
\end{equation}
\label{th-proved}
\end{theorem}
\begin{proof}
On the one hand, the Nakayama automorphism $\alpha: H \rightarrow H$ for the Frobenius
homomorphism $\es \psi$ is by Proposition (\ref{prop-F}) given
by 
\[
\alpha(a) = S^2(a) \leftharpoonup m = S^2(a \leftharpoonup m)
\]
for every $a\in H$.  On the other hand, the Nakayama automorphism $\nu$ of $H$ for the
Frobenius homomorphism $\psi \in H^*$ is by Theorem (\ref{th-C}) 
\[
\nu(a) = \es^2(m \rightharpoonup a) = m \rightharpoonup \es^2(a),
\]
for every $a \in H$. 
By Proposition (\ref{prop-G}), $\psi \circ \es = \gamma \psi b$, so by the comparison theorem
\[
\alpha(a) = b^{-1} \nu(a) b
\]
for every $a\in H$. 

Substituting the first two equations in the third yields,
\[
S^2(a) = b^{-1} \es^2(m \rightharpoonup a) b \leftharpoonup m^{-1} 
\]
which is equivalent to Eq.\ (\ref{eq:radford}) since $S^2$ fixes $b$ and $m$,
and for every group-like $a \in H$, we have $m \rightharpoonup (axa^{-1}) = a(m \rightharpoonup
x)a^{-1}$. 
\end{proof}

\begin{remark}
In \cite{KS} it was shown that a group-like element $g $ in a finite
projective Hopf algebra over a Noetherian ring $k$  
has finite order dividing the least common multiple $N$ of the local
ranks of $H$.  Since $m$ and $b$ are group-like elements in $H^*$ and $H$, respectively,
 it follows
from the general Radford formula and Eq.\ (\ref{eq:nu}) 
that the antipode $S$ and the Nakayama automorphism $\nu: H \rightarrow H$
have finite order dividing $4N$ and $2N$ respectively.


Waterhouse  sketches a different method of how to 
extend the Radford formula to a finite projective Hopf algebra
and show that $S$ has finite order \cite{Water}.
Schneider has   established
Radford's formula  by  different 
methods for $k$ = field \cite{Arg}. 
Radford's formula is generalized
to double Frobenius algebras over fields by Koppinen \cite{Kop98}.
\end{remark}

\section{When Hopf algebras are separable}

In this section we give a criterion in terms of the left norm $N$ for
when a finite projective Hopf algebra $H$ is separable.  We  first
need a proposition closely related to some results on 
when Frobenius algebras/extensions/bimodules are separable
\cite{HS,CK,NEFE}.  Let $k$ be a commutative ground ring.

\begin{proposition}
Suppose $A$ is a $P$-Frobenius algebra with system $(\psi, x_i,q_i,y_i)$.
Then $A$ is $k$-separable if and only if there is $d \in P \otimes A$
such that
\[
\sum_i x_i q_i d y_i = 1_A.
\]
\label{prop-june7}
\end{proposition}
\begin{proof}
The forward implication is proven by first letting $\sum_j a_j \otimes b_j$
be the separability element for $A$.
Next set $d := \sum_j \psi(a_j) \otimes b_j \in P \otimes H$. 
Then
\[
\sum_i x_i q_i d y_i = \sum_{i,j} x_i q_i\psi(a_j) b_j y_i =
\sum_j \sum_i x_i q_i \psi(y_i a_j)b_j = \sum_j a_j b_j = 1_A. 
\]

The reverse implication is proven by noting that $e := \sum_i  x_i \otimes q_i d y_i$
is a separability element for $A$. By hypothesis, $\mu(e) = 1$ where $\mu: A \otimes A \rightarrow
A$ is the multiplication mapping.  $e$ is in the center $(A \otimes A)^A$
of the natural $A$-bimodule $A \otimes A$ as a consequence of Proposition (\ref{prop-enumb}).
\end{proof}

Next, let $P$ be
an invertible $k$-module with inverse $Q$.
We shall say that $q \in Q$ is {\it Morita-invertible} if there is $p \in P := Q^*$ such
that $qp = 1_k$.  Note that a left inverse in this sense may differ from
a right inverse by a unit $\chi$ in $k$, since $qp = \chi pq$.
We note that if $q \in Q$ is Morita-invertible, then $Q$ and $P$ are free of rank one,
since $q' \mapsto q'p$ is epi $Q \to k$, whence an isomorphism.  
More generally, we say that $\sum_i q_i \otimes a_i \in Q \otimes A$
 is
Morita-invertible where $A$ is a $k$-algebra
 if there is $\sum_j p_j \otimes b_j \in P \otimes A$
such that $\sum_{i,j} q_ip_j a_i b_j = 1_A$.  
The next theorem generalizes  results in \cite{OS73,BFS2}.

\begin{theorem}
Suppose $H$ is a finite projective Hopf algebra with $P$-Frobenius homomorphism
$\psi$ satisfying Eq.\ (\ref{eq:li}) and left norm $N = \sum_i N_i \otimes q_i$.
Then $H$ is $k$-separable if and only if $\sum_i \epsilon(N_i)q_i$ is
Morita-invertible.
\label{th-blueskies}
\end{theorem}
\begin{proof}
We make use of the dual bases $\{ N\II \}, \{ q_i\} , \{ \es(N\I) \}$
given by Proposition (\ref{prop-par73}). 
If $H$ is $k$-separable, then by the proposition above there is $d
:= \sum_j p_j \otimes a_j \in P \otimes H$
such that 
\[
\sum_{i,(N_i)} N\II q_id \es(N\I) = 1_H.
\]
Applying $\epsilon$ we obtain
\[
\sum \epsilon (\epsilon(N\I) N\II)q_i p_j \epsilon(a_j) = \sum_i \epsilon(N_i)q_i \sum_j p_j \epsilon(a_j) = 1_k,
\]
whence $\sum_i \epsilon(N_i) q_i$ is Morita-invertible. 

Conversely, if $q := \sum_i \epsilon(N_i) q_i$ is Morita-invertible
with inverse $p \in P$ such that $qp = 1_k$,
then we let $d : = p \otimes 1_H$.  Note that
\[
\sum N\II q_id \es(N\I) = \sum_i \epsilon(N_i) q_ip 1_H = 1_H,
\]
whence $H$ is $k$-separable by Proposition (\ref{prop-june7}). 
\end{proof}

It follows directly from this theorem that a $k$-separable projective Hopf algebra
is an FH-algebra, since it is $P$-Frobenius
with $P \cong k$. 
Since a $k$-separable $H$ has a Morita invertible element, it
is $P$-Frobenius with $P \cong k$; whence the corollary below. 
   
\begin{corollary}
A separable Hopf algebra $H$ is an FH-algebra.
\end{corollary}
  
As a result, a separable Hopf algebra $H$  is unimodular \cite{KS}:
i.e. $m = \eps$. 
Then as in \cite{KS} we get the
following corollary which is an extension of the main theorem in
Etingof and Gelaki \cite{EG} to the case of rings.

\begin{corollary}
Suppose $2$ is not a zero-divisor in $k$,
and Hopf $k$-algebra $H$ is separable and coseparable.  Then $S^2 = \id_H$.  
\end{corollary}

Next we study when separable Hopf algebras are strongly separable.  
Recall that an algebra $A$ is {\it strongly separable} \cite[Kanzaki, Hattori]{Kan,Hat}
if there is $e : = \sum_j z_j \otimes w_j \in A \otimes A$ such 
that $\mu(e) = \sum_j z_j w_j = 1_A$ and for every $a \in A$, we have
$\sum_j z_j a \otimes w_j = \sum_j z_j \otimes aw_j$.  We will call such an 
$e \in A \otimes A$ a {\it Kanzaki separability element}:
one may prove that its transpose $\sum_i w_j \otimes z_j$ is an ordinary separability idempotent
\cite{Hat} (cf.\ \cite[Theorem 3.4]{KS2}).  For example, if
 $k$ is an algebraically closed field of characteristic $p$, then
$A$ is strongly separable if it is semisimple and none of its simple modules have dimension
over $k$  divisible
by $p$. 
We first  need a proposition which generalizes part of \cite[Prop.\ 4.1]{KS2}.
\begin{proposition}
Suppose $A$ is a $P$-Frobenius algebra with system $(\psi, x_i,q_i,y_i)$
such that 
\begin{equation}
u := \sum_i q_i \otimes y_i x_i
\end{equation}
is Morita-invertible.  Then $A$ is strongly separable.
\label{prop-sunray}
\end{proposition}
\begin{proof}
Suppose $\sum_j p_j \otimes a_j \in P \otimes A$ satisfies $\sum_{i,j}
q_ip_j y_i x_i a_j = 1_A$.  From this and Proposition (\ref{prop-enumb}),
we easily see that $e : = \sum_i y_i \otimes x_i q_i p_j a_j$ is a Kanzaki
separability element. 
\end{proof}

Setting  $u^{-1} : = \sum_j p_j \otimes a_j$, we can apply Proposition (\ref{prop-enumb})
to obtain a formula for the Nakayama automorphism:
\begin{equation}
\nu(a) = u a u^{-1},
\label{eq:bfs}
\end{equation}
where we make use of the usual Morita mapping $Q \otimes P \rightarrow k$. 

Recall that a Hopf algebra $H$ is {\it involutive} if $S^2 = \id_H$.
The next theorem contains a result of Larson \cite{Larson71} as a special case. 

\begin{theorem}
Suppose $H$ is a finite projective, separable, involutive Hopf algebra.
Then $H$ is strongly separable.
\end{theorem}
\begin{proof}
If $(\psi, N\II,q_i, \es(N\I))$ is the $P$-Frobenius system for $H$
given by 
Proposition (\ref{prop-par73}),  we note here that $\es = S$, so that
the $u$-element of Proposition (\ref{prop-sunray}), 
\[
u := \sum_i q_i \otimes \sum_{(N_i)} S(N\I)N\II = \sum_i q_i \epsilon(N_i) \otimes 1_H
\]
is Morita-invertible by Theorem (\ref{th-blueskies}). 
\end{proof}

\section{Hopf Subalgebras}

Throughout this section, $k$ is a commutative ring and 
we consider  a finite projective Hopf algebra $H$ with Hopf subalgebra $K$ which is also
finite projective as a $k$-module.
We will show that the functors of induction and co-induction from the category
${\mathcal M}_K$ of $K$-modules to ${\mathcal M}_H$ are naturally isomorphic
 up to a Morita auto-equivalence of ${\mathcal M}_K$ determined by a relative
Nakayama automorphism and a relative Picard group element. This section
generalizes results in \cite{OS73,S,KS}. 

Let $R$ be an arbitrary ring, $\beta: R \rightarrow R$ a ring automorphism,
and $M_R$ a module over $R$.  The $\beta$-twisted module $M_{\beta}$
is defined by $m \cdot r := m \beta(r)$, clearly another $R$-module.
If $\beta$ is an inner automorphism, is easy to check that $M_R \cong M_{\beta}$ and
 $M \otimes_R R_{\beta} \cong M_{\beta}$.
Then  
the bimodule ${}_RR_{\beta}$ induces a Morita auto-equivalence
of ${\mathcal M}_R$ via tensoring. 

\begin{lemma}
If $A$ is a $P$-Frobenius $k$-algebra with Frobenius homomorphism $\phi$
and corresponding Nakayama automorphism $\nu$, then we have the following
bimodule isomorphisms: 
\begin{equation}
{}_A A_A \cong {}_A\Hom(A,P)_{\nu} \cong {}_{\nu^{-1}}\Hom(A,P)
\label{eq:Yamamoto}
\end{equation}
\label{lemma-Yama}
\end{lemma}
\begin{proof}
Since $a\phi = \phi \nu(a)$ in $A^*$ for every $a\in A$,
it follows that the Frobenius isomorphisms $a \mapsto a \phi$ and
$a \mapsto \phi a$ induce the first and second isomorphisms above
(between $A$ and $\Hom(A,P)$). 
\end{proof}

  As a straightforward extension of Definition (\ref{def-PF}), we define
$P$-Frobenius extension $A/S$, where $P$ is an invertible $S$-bimodule
(and $-\otimes_S P$ defines a Morita auto-equivalence of ${\mathcal M}_S$ \cite{Bass}). 

\begin{definition}
Suppose $S$ is a subring of ring $A$ and $P$ is an invertible $S$-bimodule.
We say $A$ is a {\it $P$-Frobenius extension of $S$, or $A/S$ is 
a Frobenius extension of the third kind,} if
\begin{enumerate}
\item $A_S$ is a finite projective module;
\item ${}_SA_A \cong {}_S\Hom_S(A_S, P_S)_A$
\end{enumerate}
\end{definition} 

A $P$-Frobenius extension has a symmetric definition, a Frobenius system like in 
Section~2, a Nakayama automorphism defined on the centralizer subalgebra $C_S(A)$ of $A$ \cite{Par73},
and a comparison theorem, which we will not need here.  As a straightforward
consequence of a theorem by Morita \cite{Mor65,Mor67}, we state without
proof (cf.\ \cite{FMS}): 
\begin{proposition}
$A$ is $P$-Frobenius extension of $S$ if and only if there is a natural isomorphism
of right $A$-modules, 
\begin{equation}
M \otimes_S A \cong \Hom_S(A_S,M \otimes_S P_S)
\end{equation}
for every module $M \in {\mathcal M}_S$. 
\end{proposition}
 
This equivalent condition for a $P$-Frobenius extension states in other words
that the functors of induction and co-induction from ${\mathcal M}_S$
into ${\mathcal M}_A$ form a commutative triangle with the Morita auto-equivalence
of ${\mathcal M}_S$ induced by $-\otimes_S P$. 

Suppose a Frobenius algebra pair forms a projective ring
extension such that the Nakayama automorphism of the overalgebra preserves the
subalgebra.  We now obtain a theorem that states that such a pair forms a 
certain $P$-Frobenius extension. 

\begin{theorem}
Suppose $A$ is a $P$-Frobenius algebra, $B$ is a $P'$-Frobenius algebra,
and $B$ is subalgebra of $A$ such that $A_B$ is a finite projective module,
and a Nakayama automorphism $\nu_A$ of $A$ sends $B$ into $B$: $\nu_A(B) = B$.
Let $\nu_B$ denote a Nakayama automorphism of $B$. Then 
$A$ is a $W$-Frobenius extension of $B$, where 
\begin{equation}
W = {}_{\beta}B \otimes Q' \otimes P, 
\label{eq:doubleu}
\end{equation}
$Q' = P^{'*}$ and $\beta$ is the relative Nakayama automorphism given by 
\begin{equation}
\beta = \nu_B \circ \nu_A^{-1}.
\label{eq:beta}
\end{equation}
\label{th-parmor}
\end{theorem}
\begin{proof}
Since $A_B$ is assumed finite projective, we need only show that ${}_BA_A \cong
{}_B\Hom_B(A_B, W_B)$.  We compute using the hom-tensor adjointness relation
and two applications of Lemma (\ref{lemma-Yama}):
\begin{eqnarray*}
{}_BA_A & \cong &  {}_{\nu_A^{-1}}\Hom(A,P)_A \\
& \cong & \Hom_k(A \otimes_B B_{\nu_A^{-1}}, k)_A \otimes P \\
& \cong & {}_{\nu_A^{-1}}\Hom_B(A_B, B^*_B)_A \otimes P \\
& \cong & {}_{\nu_A^{-1}}\Hom_B(A_B, {}_{\nu_B}B_B \otimes Q')_A \otimes P \\
& \cong & {}_B\Hom_B(A_B,\, {}_{\nu_B \circ \nu_A^{-1}} B_B \otimes Q' \otimes P)_A \qed
\end{eqnarray*}
\renewcommand{\qed}{}\end{proof}

Let $K  \subseteq H$ be a  pair of finite projective 
Hopf $k$-algebras where $K$ is a Hopf subalgebra
of $H$ (i.e., $K$ is a pure
$k$-submodule of $H$,  $\Delta(K) \subseteq K \otimes K$ and $S(K) = K$) in the next corollary.
Let $P(K)^*$, $P(H)^*$ be the $k$-module of integrals $\int^{\ell}_K$, $\int^{\ell}_H$,
respectively,  $\nu_H$, $\nu_K$ be the respective Nakayama automorphisms
and $m_H$, $m_K$ be the respective left modular functions.  

\begin{corollary}
If $K \subseteq H$ is a finite projective Hopf subalgebra pair,  then $H/K$
is a $P$-Frobenius extension where  
\begin{equation}
P = {}_{\beta}K \otimes P(K)^* \otimes P(H)
\label{eq:pee}
\end{equation}
and 
\begin{equation}
\beta = \nu_K \circ \nu_H^{-1}.
\label{eq:bee}
\end{equation}
\end{corollary}
\begin{proof}
The natural module $H_K$ is finite projective as a corollary of the Nichols-Zoeller
Freeness theorem \cite[Prop.\ 5.3]{KS}. 
Furthermore, the Nakayama automorphism
$\nu_H^{\pm 1}(a) = m_H^{\pm 1} \rightharpoonup S^{\mp 2}(a)$ for 
every $a \in H$
by Eq.\ (\ref{eq:nu}), whence $\nu_H(K)=  K$. Thus the hypotheses of Theorem (\ref{th-parmor})
are satisfied.
\end{proof}

It follows from the formulas for $\nu_H$ and $\nu_K$ in Eq.\ (\ref{eq:nu})
that for every $x \in K$, 
\begin{eqnarray}
\beta(x) & = &  m_K \rightharpoonup \es^2(m_H^{-1} \rightharpoonup S^2(x)) \nonumber \\
& = & (m_K * m_H^{-1}) \rightharpoonup x
\label{eq:relmod}
\end{eqnarray}
(cf. \cite{FMS}). 



 \begin{remark}
Kasch  makes a study in \cite{Kas2} of  the relative homological
algebra of Frobenius extensions.  One can extend this study 
 to a Frobenius
extension $A/S$ of the third kind by taking into account some Morita theory.
  For example, one may show by these means
that under the (rather common) additional assumption that $S$ is $S$-bimodule isomorphic
to a direct summand in $A$,  the flat dimension of any
 $S$-module is equal to both
the flat dimension of its induced $A$-module  and of its
co-induced $A$-module. In  \cite{Par72} the study in \cite{Kas2} 
is extended to a  cohomology
theory for FH-algebras, showing that 
these  have a complete cohomology with
cup product, a generalized Tate duality
under a certain cocommutativity condition,
 and a generalized Hochschild-Serre spectral sequence.  
\end{remark}

\section{Embedding $H$ into an FH-algebra}

In this section we show that a finite projective Hopf algebra $H$ is a Hopf subalgebra
of an FH-algebra in two ways.  We first show that $H$ is a Hopf subalgebra of $D(H)$.
We let $k$ be a commutative ring.
The {\it quantum double} $D(H)$ of a finite dimensional Hopf
algebra, due to Drinfel'd
\cite{Drin}, is readily extended to a  finite projective 
Hopf algebra $H$ over $k$:  
 at the level of coalgebras it is given by  
\[
D(H) := H^{*\, {\rm cop}} \otimes_k H ,
\]
where $H^{*\, {\rm cop}}$ is the co-opposite of $H^*$, the coproduct
being $\Delta^{\rm op}$. 
 The multiplication
on $D(H)$ is described in two equivalent ways
 as follows  \cite[Lemma 10.3.11]{Mont}. In terms of  the notation
 $gx$ replacing $g \otimes x$ for every $g \in H^*, x \in H$, both $H$ and $H^{*\, {\rm cop}}$ are subalgebras of
$D(H)$, and for each $g \in H^*$ and $x \in H$,
\begin{equation}
xg := \sum(x\1 g S^{-1}x\3)x\2 
= \sum g\2(S^{-1}g\1 \rightharpoonup x \leftharpoonup g\3).
\end{equation}
The algebra $D(H)$  is a Hopf algebra 
 with antipode $S'(gx) :=  Sx S^{-1}g$, the proof of this proceeding as in 
\cite{CK2}.

\begin{theorem}
If $H$ is a finite projective Hopf algebra, then $D(H)$ is an FH-algebra.
\end{theorem}
\begin{proof}
It is enough to show that $\int^{\ell}_{D(H)^*} \cong k$. 
As an algebra, $D(H)^* \cong H^{\rm op} \otimes H^*$, the tensor product algebra
of $H^*$ and the opposite algebra of $H$. Now $H$
is $P$-Frobenius algebra if and only if $H^{\rm op}$ is, since they have the same
Frobenius system with a change of order in the dual base.  By Proposition (\ref{prop-fixed}),
 $H^*$ is a $P^*$-Frobenius algebra. 
It follows from Lemma (\ref{lemma-june1}) that $D(H)^*$ is a Frobenius algebra, since
$P \otimes P^* \cong k$. Now the $k$-space of integrals of an
 augmented Frobenius algebra is free of rank one, which proves our
theorem. 
\end{proof}

Next we show that $H$ has a ring extension to an FH-algebra $H \otimes_k K$.
This will follow right away from the construction of a ring extension
$k \subset K$ where $K$ has trivial Picard group. 
We continue with $k$ as a commutative ring, and let $M$ be the set of all maximal
ideals in $k$.  Choose finite subsets $M_{\alpha} \subset M$, $\alpha \in I$
 such that $ \cup_{\alpha \in I}
M_{\alpha} = M$ and the subsets $M_{\alpha}$ are linearly ordered with respect to inclusion:
in other words, for any two indices $\alpha, \beta \in I$ either  $M_{\alpha} \subset M_{\beta}$
or $M_{\beta} \subset M_{\alpha}$.

Let $m_{\alpha_1}, \ldots, m_{\alpha_n}$ be all the elements  of $M_{\alpha}$,
i.e. maximal ideals in $k$. Then the
set \[
K_{\alpha} = k_{m_{\alpha_1}}\oplus \cdots \oplus k_{m_{\alpha_n}}
\]
 is a semilocal ring and has
trivial Picard group: $Pic(K_{\alpha}) = 0$. 
For any pair $M_{\alpha} \subset M_{\beta}$, we have the canonical projection $\pi_{\alpha \beta}:
K_{\beta} \rightarrow K_{\alpha}$ and we may consider the inverse limit ring 
\begin{equation}
K := \lim_{\leftarrow} (K_{\alpha},\pi_{\alpha \beta}) .
\label{eq:aas}
\end{equation}
Furthermore, for any $\alpha \in I$ we have the canonical homomorphism $f_{\alpha}: k \rightarrow K_{\alpha}$,
which is the direct sum of the corresponding localization homomorphisms.
The following diagram is clearly commutative:

$$\begin{diagram}
& & k & & \\
& \SW^{f_{\beta}} & & \SE^{f_{\alpha}}  & \\
K_{\beta} & & \rTo^{\pi_{\alpha \beta}}& & K_{\alpha} 
\end{diagram}$$
\vspace{2mm}
From universality we obtain a homomorphism $f: k \rightarrow K$.

\begin{lemma}
$f$ is a monomorphism.
\end{lemma}
\begin{proof}
Let $f_m$ be the localization homomorphism $f_m: k \rightarrow k_m$.
Then it follows easily that $\ker f = \cap_{m \in M} \ker f_m = 0$.
\end{proof}

Now let $\pi_{\alpha}: K \rightarrow K_{\alpha}$ be the canonical epi.
Since the diagram

$$\begin{diagram}
& & K & & \\
& \SW^{\pi_{\beta}} & & \SE^{\pi_{\alpha}}  & \\
K_{\beta} & & \rTo^{\pi_{\alpha \beta}} & & K_{\alpha} 
\end{diagram}$$
\vspace{2mm}
is commutative, the following diagram is commutative as well:

$$\begin{diagram}
& & Pic(K) & & \\
& \SW^{ Pic(\pi_{\beta})} & & \SE^{Pic(\pi_{\alpha})} & \\
Pic(K_{\beta}) & & \rTo^{Pic(\pi_{\alpha \beta})} & & Pic(K_{\alpha}) 
\end{diagram}$$
\vspace{2mm}
Again from universality we obtain a homomorphism
\[
\Phi:\ Pic(K) \longrightarrow \lim_{\leftarrow}(Pic(K_{\alpha}), Pic(\pi_{\alpha \beta}))
\]
\begin{theorem}
$\Phi$ is injective.
\end{theorem}
\begin{proof}
We need the following result proved in \cite{BS}:
\begin{theorem}
Suppose $I$ is some linearly
ordered  set  and 
 for each ordered $\alpha, \beta \in I$, $A_{\alpha}$ is a commutative ring and 
there is an epimorphism $\psi_{\alpha \beta}$ such that
the restriction to the group of units $\psi_{\alpha \beta}: U(A_{\beta}) \rightarrow U(A_{\alpha})$
is a surjection. If 
$$A = \lim_{\leftarrow}( A_{\alpha}, \psi_{\alpha \beta}),$$ then the induced
map 
$$Pic(A) \rightarrow \lim_{\leftarrow} (Pic(A_{\alpha}), Pic(\psi_{\alpha \beta}))$$
is injective.
\end{theorem}
The hypotheses of this proposition are fulfilled by the mappings
$\pi_{\alpha \beta}: K_{\beta} \rightarrow K_{\alpha},$
whence $\Phi$ is injective. 
\end{proof}

The next corollary follows from recalling that $Pic(K_{\alpha}) = 0$. 
\begin{corollary}
Given a commutative ring $k$
and $K$ defined in Eq.\ (\ref{eq:aas}), $k \subset K$ is a ring
extension with  $Pic(K) = 0$.
\end{corollary}

\end{document}